\documentclass{article}
\def\uu{\bigsqcup}

\def \Z {{\mathbf {Z}}}
\def \T {{\mathbf {T}}}
\def \N {{\mathbf {N}}}

\title{ Кратности сверток спектров автоморфизмов}
\author{ Valery V. Ryzhikov}
\date{}

\usepackage[T1]{fontenc} 
\usepackage[cp1251]{inputenc}  
\usepackage[russian]{babel}
\begin{document}

\begin{center}   
{ \Large \bf
Spectral properties of dynamical tensor powers, \\ and tensor factorizations of simple Lebesgue spectrum 
    
\vspace{5mm}

Valery V. Ryzhikov}
\end{center} 

\vspace{4mm} 
{\large For every $n>0$ there is a unitary operator $U$ such that the shift in $l_2(\Z)$  with  simple Lebesgue spectrum is isomorphic to the tensor product $U\otimes U^2\otimes\dots\otimes U^{2^n}.$ There is an ergodic automorphism $T$ with its  symmetric tensor power $T^{\odot n}$ of   simple spectrum, and  $T^{\odot(n+1)}$ of    absolutely continuous spectrum.}



 

\Large
\section{Introduction}
A unitary operator $U$ with  simple spectrum (in other words, with spectrum of multiplicity 1) is equivalent to the operator of multiplication by a variable $z\in \mathbf C$, $|z| = 1$, acting in space
$L_2(\T, \sigma)$, where $\sigma$ is a Borel measure on the unit circle $\T = \{t : |t| = 1\}$. Such a measure and any equivalent to it
is called the spectral measure of the operator $U$. Direct sums (finite and countable) of such operators give a spectral representation for all possible unitary operators $V$ in the separable Hilbert space.

 The unitary operator is  defined by
the measure of the maximum spectral type $\sigma_{max}$ and the function of spectral multiplicity ${\cal M}:\T\to \N\cup \{\infty\}$, measurable w.r. to $\sigma_{max}$.

One of the main questions in the spectral theory of actions with an invariant measure is: \it which pairs $(\sigma_{max},{\cal M})$ correspond to the spectral representation of an ergodic automorphism $T$? \rm

Is there an ergodic realization of the pair $(\sigma_{max}, {\cal M})$
provided that
 $\sigma_{max}$ is absolutely continuous, but the multiplicity function is bounded? This is unknown.

The special case when ${\cal M}\equiv 1$ and $\sigma_{max}=\lambda$ ( the Lebesgue measure) is known as  Banach's problem: in other words, is there an automorphism with  simple  Lebesgue spectrum?

The Banach problem consists of two problems. In Ulam's book he is talking about an automorphism of a space with infinite measure (isomorphic to a line with the Lebesgue measure). In the case of an automorphism of a probability space, we consider its action on the space of functions from $L_2$ with zero mean.

Candidates for the role of automorphisms with  simple Lebesgue spectrum are automorphisms of rank 1 (see \cite{B}, \cite{R14}). Jean-Paul Thouvenot, however, suggested me that the examples are outside the class of rank 1 transformations.
Without completely abandoning transformations of rank 1, we can still assume that a suitable example exists among products of the form $S\otimes T$ and even among $T\odot T$, where $S,T\in$ Rank-1. Such $S,T$ must have  singular spectra. Here, of course, only automorphisms of space with infinite measure are considered (for them, nonzero constants are not included in $L_2$, and the product $S\otimes T$ has no coordinate factors). But even in the case of the  probability space, the ergodic realization for the product $S\otimes T$ of  simple spectrum with an absolutely continuous component would be a new important result for the spectral theory of dynamical systems.

Thus we come to the question of
compatibility of the following properties  of ergodic automorphisms $S,T$  (the case $T=S$ is also of interest):

\bf (a.c.) \rm the spectrum of $S\otimes T$ is not singular;

\bf (f.m.) \rm the spectrum $S\otimes T$ is of finite multiplicity.

It is natural to first ask the question: \it is there a pair of unitary operators $S,T$ with continuous spectrum that satisfy
\bf (a.c.) + (f.m.) \rm ?

\vspace{3mm}
The answer is yes even for $T=S^2$.

\vspace{3mm}
\bf Theorem 1. \it For every $n>0$ there is a unitary operator $U$ such that the unitary operator with  simple Lebesgue spectrum is isomorphic to the tensor product
$$U\otimes U^2\otimes\dots\otimes U^{2^n}.$$ \rm

\vspace{3mm}
Perhaps the reader will prefer to postpone reading this note and find the tensor factorization of the standard unitary shift in $l_2(\Z)$ in his own way, and then compare it with our approach. 

\vspace{3mm}
Theorem 1 does not shed light on the Banach problem, since the spectra of operators of dynamic origin have their own specifics, but it can add enthusiasm to the study  spectra of automorphism tensor products.

As an example of some  unusual behavior of the spectra of symmetric tensor powers, we present the following statement.

\vspace{3mm}
\bf Theorem 2. \it For $n>1$ there is an ergodic automorphism $T$ such that the symmetric tensor power $T^{\odot n}$ has  simple spectrum, and the power $T^{\odot (n+1 )}$ has  Lebesgue spectrum (the case $n=1$ is considered in \cite{R24}). \rm

\vspace{3mm}
We will describe the constructions of such automorphisms and communicate the idea   proof of Theorem 2, restricting ourselves to the case $n=2$.

\section{ Constructions}
We fix the natural number $h_1$, the sequence $r_j$ (the number of columns into which the tower of stage $j$ is cut) and the sequence of integer vectors (spacer parameters)
$$ \bar s_j=(s_j(1), s_j(2),\dots, s_j(r_j-1),s_j(r_j)).$$

At step $j$, a system of disjoint half-intervals is defined
$$E_j, TE_j, T^2E_j,\dots, T^{h_j-1}E_j,$$
and on half-intervals $E_j, TE_j, \dots, T^{h_j-2}E_j$
the transformation $T$ acts by parallel translation. Such a set of half-intervals is called a tower of stage $j$; their union is denoted by $X_j$ and is also called a tower.

Let us represent $E_j$ as a disjoint union of $r_j$ half-intervals
$$E_j^1,E_j^2E_j^3,\dots E_j^{r_j}$$ of the same length.
For each $i=1,2,\dots, r_j$ we define the column $X_{i,j}$ as the union of intervals
$$E_j^i, TE_j^i ,T^2 E_j^i,\dots, T^{h_j-1}E_j^i.$$

To each column $X_{i,j}$ we add $s_j(i)$ of disjoint half-intervals of the same measure as $E_j^i$, obtaining a set
$$E_j^i, TE_j^i, T^2 E_j^i,\dots, T^{h_j-1}E_j^i, T^{h_j}E_j^i, T^{h_j+1}E_j^i , \dots, T^{h_j+s_j(i)-1}E_j^i$$
(all these sets do not intersect).
Denoting $E_{j+1}= E^1_j$, for $i<r_j$ we set
$$T^{h_j+s_j(i)}E_j^i = E_j^{i+1}.$$
 The set of superstructured columns is from now on considered as a tower of stage $j+1$, consisting of half-intervals
$$E_{j+1}, TE_{j+1}, T^2 E_{j+1},\dots, T^{h_{j+1}-1}E_{j+1},$$
Where
 $$ h_{j+1}+1 =(h_j+1)r_j +\sum_{i=1}^{r_j}s_j(i).$$

The partial definition of the transformation $T$ at step $j$ is preserved in all subsequent steps. As a result, an invertible transformation $T:X\to X$ is defined on the space $X=\cup_j X_j$, preserving the standard Lebesgue measure on $X$.

Transformation $T$ and its induced unitary operator
$T$, $Tf(x)=f(Tx)$, are denoted identically in the article.
An automorphism of rank one is ergodic and has  simple  spectrum. It is known that the tower indicator $X_1$ is a cyclic vector for the operator $T$ with continuous spectrum, which is certainly true if the measure $X$ is infinite.

\vspace{2mm}
\bf An example of  a construction $\bf T$  for which 
 $\bf T^{\otimes 2}$ has spectral multiplicity 2, and 
spectrum of $\bf T^{\otimes 3}$ is absolutely continuous.

\rm
 Let the sequence $J_k$ be defined as follows:

$J_1=1,\ J_{k+1}= J_k +(k+1)^8.$

\vspace{2mm}
\bf The parameters of $T$:\rm

$r_{j}:=(k+1)^2 $, for $j(k)\leq j<j(k+1)$;

$s_j(1):=2h_j-1$; \ $s_j(2):=2h_j$;

 $s_j(i):=2^ih_j$ for $2 < i \leq r_j,$

where $h_1=2$, $h_{j+1}= r_jh_j +\sum_{i=1}{r_j} s_j(i)$.

 The value $r_j=(k+1)^2$ is repeated by   many times, so for
 the  sums $$F_k=\sum_{J_k<j< J_{k+1}} {\bf 1}_{T^{h_j}E_1}\otimes {\bf 1}_{T^{h_j}E_1} $$ the 
functions $F_k /\|F_k\|$ were asymptotically close to the indicator of 
the set $$\uu_{0\leq p,q\leq 1}{T^pE_1}\times {T^qE_1}.$$
 This is established by slightly modifying the methods of \cite{R24}. Thanks to such  approximations, we obtain the following.

\vspace{3mm}
\bf Lemma. \it Let $f={\bf 1}_{E_1}\times {\bf 1}_{E_1}$.
Then for any $m,n\in \Z$ any closed and $T\otimes T$-invariant space  contains, together with the vector $(T^m\otimes T^n)f$ all 
vectors $$(T^{m+1}\otimes T^n)f + (T^m\otimes T^{n+1})f.$$\rm

\vspace{3mm}
From Lemma it now easily follows that $T\otimes T$-invariant space containing the vectors $f$ and $(T\otimes I)f$,
contains all vectors of the form $(T^m\otimes T^n)f$. This means that $T\otimes T$ has a homogeneous spectrum of multiplicity 2. Consequently,   $T\odot T$ has  simple spectrum, since
$$T\otimes T\cong (T\odot T) \,\bigoplus\, (T\odot T).$$

If $\hat\sigma(n)=\mu(E_1\cap T^nE_1)$, then for the coefficients
$$c(n):=\hat\sigma(n)^3=\hat{\sigma^{\ast 3}}(n)$$
the convergence of the series $\sum_n c(n)^2$ can be  directly verified. Consequently, the spectrum $T^{\odot 3}$ is absolutely continuous, since the spectral type of the operator $T^{\odot 3}$ is the convolution power $\sigma^{\ast 3}$ of the spectral measure $\sigma$ for  the operator $ T$. If the convergence of the indicated series is fast enough, then $\sigma^{\ast 3}$ is guaranteed to be equivalent to the Lebesgue measure.

\vspace{2mm}
Numerous results on spectra of measure-preserving actions are described in  \cite{KL}. To these the author can, for example, add constructive solutions to the Kolmogorov problem on the group property of the spectrum and the Rokhlin problem on the homogeneous spectrum in the class of mixing systems; solution of the Bergelson spectral question on the compatibility of rigid and mixing sequences; the answer to question
of Oseledets on projections of a tensor square (spectral) measure.

\section {Proof of Theorem 1} 
\bf Case $\bf n=2$. \rm On $[0, 1]$ consider
two sets $A,B$:

$$A=\left\{\sum_{i=1}^\infty \frac {a_i} {2^{2i}} \, :\, a_i\in\{0,1\}\right\},$$

$$B=2A=\left\{\sum_{i=1}^\infty \frac {2a_i} {2^{2i}} \, :\, a_i\in\{0,1\}\right\}. $$

 Now we define  two measures with supports $A$ and $B$, respectively.
The first $\sigma_1$ is the  Bernoulli measure of type
$(\frac 1 2, \frac 1 2, 0,0)$, the second measure $\sigma_2$ is also Bernoulli of type $(\frac 1 2, 0, \frac 1 2, 0)$.

(i) Mapping $C_2: A\times B\to [0,1]$,
$$ C_2(a,b)=a+b, $$
we get  an isomorphism of the spaces $([0,1]^2, \sigma_1\times \sigma_2)$
and $([0,1],\lambda)$ ($\lambda$ is the Lebesgue measure). This is easy to check (the standard method: the cylinder goes into the interval while maintaining the measure).

Let $U$ be a multiplication operator acting in the space $L_2([0,1], \sigma_1)$ as follows:
$$ Uf(x)=e^{2\pi i x}f(x).$$

Then the action of the operator $V=U\otimes U^2$ on $([0,1]^2, \sigma_1\times \sigma_1)$ is described by the formula
$$ VF(x,x')=e^{2\pi i (x+2x')}F(x,x').$$
But for almost every $y\in [0,1]$ there are unique $x,x'\in A$ such that $x+2x'=y$, so by (i) the operator $V$ is isomorphic
to the operator of multiplication by $z=e^{2\pi i y}$ in the space $L_2([0,1].\lambda)$, which is what needed to be proved.

\vspace{2mm}
\bf Case $\bf n=3$.  \rm
Let's define
$$A=\left\{\sum_{i=1}^\infty \frac {x_i} {2^{3i}} \, :\, x_i\in\{0,1\}\right\};$$

We note that for almost every $y$ there are unique $x,x',x''\in A$
such that $y=x+2x'+4x''$, and the mapping
$$(x,x',x'')\to x+2x'+4x''$$
gives  an isomorphism of spaces
$([0,1]^3, \sigma_1^{\otimes 3})$ and
and $([0,1],\lambda)$, where
$\sigma_1$ -- Bernoulli measure of type $(\frac 1 2, \frac 1 2, 0.0, 0, 0, 0.0)$,

Let $U$ now denote the multiplication operator in the space $L_2([0,1], \sigma_1)$:
$$ Uf(x)=e^{2\pi i x}f(x).$$

Then the operator $V_3=U\otimes U^2\otimes U^4$ is 
the multiplication   by  $e^{2\pi i (x +2x'+4x'')}$ in the space   $([0,1]^3, \sigma_1^{\otimes 3})$, so it has  simple Lebesgue spectrum.

 \vspace{4mm}
Similar reasoning for $\bf n>3$. On the set $A$,
$$A=\left\{\sum_{i=1}^\infty \frac {x_i} {2^{ni}} \, :\, x_i\in\{0,1\}\right\},$$ 
 we consider halph-halph Bernoulli measure $\sigma_1$, then we see
that the multiplication   by  $e^{2\pi i (x +2x^1+\dots + 2^nx^n)}$
is an operator on $L_2(\sigma_1^{\otimes (n+1)})$ with simple Lebesgue spectrum.
Thus,  we get the operator $U$,  $Uf(x)=e^{2\pi i x}f(x)$, that acts in $L_2(\sigma_1)$, and  the product $U\otimes U^2\otimes\dots\otimes U^{2^n}$ is of simple Lebesgue spectrum.

Theorem 1 is proved.

 \vspace{4mm}
\bf Remarks. \rm 
We know how the spectral multiplicities of powers $V^n$ grow:
they are equal to $n$. We hope the reader will agree that it is interesting to know how the spectral multiplicities of the powers of $U^n$ grow?  They cannot grow asymptotically linearly, since in a tensor product we would see a  superlinear  growth. Similar question
on the structure of the spectrum of the operator 
$$Gauss(U)^n=\exp(U^n)= \bigoplus_{d=0}^\infty
 \left(U^{\odot d}\right)^n$$
is interesting as well. 

In \cite{R22} one of  Oseledets' problems is solved: there exists a (spectral) measure $\sigma$ such that for some dense set of directions  the corresponding projections of $\sigma\times\sigma$ (on the horizontal line)  are absolutely continuous. At the same time for  other directions that  are dense too  the corresponding projections along these  directions  are singular measures.

 Let some projection of the measure $\sigma\times\sigma$ be absolutely continuous, and all other projections that are sufficiently close  in direction give singular measures.  Osledets conjectured the impossibility of such a situation.  
Our example of $sigma_1$, due to its self-similarity  does not provide a counterexample,  but perhaps a little similar to it.

Is there an operator  $U$ with  $U^{\odot n}$ of   simple    absolutely continuous spectrum?

\newpage

\begin{center}   
{ \Large \bf
Спектральные свойства  тензорных степеней преобразования \\ 
и тензорные факторизации простого лебеговского спектра 
\vspace{5mm}

В.В. Рыжиков}
\end{center} 

\vspace{5mm}
\bf Введение.  \rm Унитарный оператор $U$ c простым спектром (т.е. спектральная  кратность равна   1, иначе говоря, у оператора есть циклический вектор) эквивалентен оператору умножения на переменную $z\in \mathbf C$, $|z| = 1$, действующему в пространстве
$L_2(\T, \sigma)$, где  $\sigma$ – борелевская мера на единичной окружности $\T = \{t : |t| = 1\}$. Такую меру и всякую эквивалентную ей 
называют спектральной мерой оператора $U$. Прямые суммы (конечные и счетные) таких операторов дают спектральное представление для всевозможных унитарных операторов $V$ в сепарабельном гильбертовом пространстве.

  Унитарный оператор полностью определяется 
мерой максимального спектрального типа $\sigma_{max}$ и функцией
кратности спектра  ${\cal M}:\T\to \N\cup \{\infty\}$, измеримой отностительно $\sigma_{max}$.

Один из главных вопросов спектральной теории действий с инвариантной мерой звучит так:  \it какие пары $(\sigma_{max},{\cal M})$ отвечают спектральному предствлению эргодического  автоморфизма $T$? \rm

Есть ли эргодическая реализация пары $(\sigma_{max}, {\cal M})$
при условии, что 
$\sigma_{max}$  абсолютно непрерывна, а функция  кратности ограничена? Это неизвестно.

Специальный случай, когда ${\cal M}\equiv 1$ и $\sigma_{max}=\lambda$  -- мера Лебега, известен как проблема Банаха: иначе говоря, существует ли автоморфизм с простым (однократным) лебеговским спектром?

Проблема Банаха состоит из  двух задач. В  в    книге Улама  речь идет об автоморфизме пространства с бесконечной мерой (изоморфному прямой с мерой Лебега).  В случае  автоморфизма вероятностного пространства,  рассматривается его  действие  на пространстве функций из $L_2$ с нулевым средним.  

Кандидатами на роль автоморфизмов с однократным лебеговским спектром   являются  автоморфизмы ранга 1 (см.  \cite{B}, \cite{R14}).  Жан-Поль Тувено, однако,  высказывал предположение о том, что примеры находятся вне класса преобразований ранга 1.   
Не расставаясь полностью с преобразованиями ранга 1, можно предположить, что подходящий пример есть   среди произведений  вида $S\otimes T$ и даже среди $T\odot T$, где $S,T\in Rank-1$. Такие $S,T$    должны  обладать  сингулярным спектром. Здесь, конечно, рассматриваются только автоморфизмы пространства с бесконечной мерой (для них ненулевые  константы не входят в $L_2$,   координатные факторы у произведния $S\otimes T$ отсутствуют). Но и в случае вероятностного пространства эргодическая реализация для произведения $S\otimes T$  простого спектра  с абсолютно непрерывной компонентой  была бы  новым важным  результатом для  спектральной теории динамических систем.

Таким образом, мы подошли  к вопросу  о
совместимости следующих  свойств  пары $S,T$ эргодических автоморфизмов  (случай $T=S$ также предствавляет интерес): 

\bf (a.c.) \rm спектр $S\otimes T$ имеет абсолютно непрерывную компоненту,

\bf (f.m.)  \rm спектр $S\otimes T$ конечнократный.

Естественно предварительно задаться  вопросом: найдется ли пара  унитарных операторов $S,T$ с непрерывным спектром, удовлетворяющих
\bf (a.c.) +  (f.m.)\rm ?

Ответ положительный  даже при $T=S^2$.  

\vspace{3mm}
\bf Теорема 1. \it Для всякого $n>0$ найдется унитарный оператор $U$ такой, что унитарный оператор с простым лебеговским спектром  изоморфен тензорному произведению
$$U\otimes U^2\otimes\dots\otimes U^{2^n}.$$   \rm

 Быть может,  читатель предпочтет  отложить чтение этой заметки и найти тензорную факторизацию стандартного  унитарного сдвига в $l_2(\Z)$ своим  способом,  а потом  сравнить его с нашим подходом. Поэтому
доказательство  теоремы  изложим  чуть позже.

\vspace{3mm}
Теорема 1 не проливает свет на проблему Банаха, так как спектры операторов динамического происхождения имеют свою специфику, но может добавить энтузиазма в    изучении  спектров тензорных произведений автоморфизмов.

В качестве примера непривычного  поведения спектров симметрических тензорных степеней  приведем следующее утверждение. 

\vspace{3mm}
\bf Теорема 2. \it Для $n>1$ найдется эргодический автоморфизм $T$ такой, что симметрическая тензорная степень $T^{\odot n}$ имеет простой спектр, а степень $T^{\odot (n+1)}$ имеет лебеговский спектр (случай $n=1$ рассмотрен в \cite{R24}).  \rm

\vspace{3mm}
Ограничившись случаем $n=2$, опишем соответстующую конструкцию и  сообщим схему доказательства.

\vspace{2mm}
\bf   Конструкции автоморфизмов ранга один. \rm   
По традиции напоминаем нужное нам определение класса конструкций. Фиксируем натуральное  $h_1$, последовательность  $r_j$ (число колонн, на которые  разрезается башня этапа $j$) и  последовательность целочисленных векторов (параметров надстроек)   
$$ \bar s_j=(s_j(1), s_j(2),\dots, s_j(r_j-1),s_j(r_j)).$$  

На шаге $j$  определена система   непересекающихся полуинтервалов 
$$E_j, TE_j, T^2E_j,\dots, T^{h_j-1}E_j,$$
а на полуинтервалах $E_j, TE_j, \dots, T^{h_j-2}E_j$
пребразование $T$ действует  параллельным переносом. Такой набор   полуинтервалов  называется башней этапа $j$, их объединение обозначается через $X_j$ и также называется башней.

Представим   $E_j$ как дизъюнктное объединение  $r_j$ полуинтервалов 
$$E_j^1,E_j^2E_j^3,\dots E_j^{r_j}$$ одинаковой длины.  
Для  каждого $i=1,2,\dots, r_j$ определим  колонну $X_{i,j}$ как объединение интервалов  
$$E_j^i, TE_j^i ,T^2 E_j^i,\dots, T^{h_j-1}E_j^i.$$

К каждой  колонне $X_{i,j}$ добавим  $s_j(i)$  непересекающихся полуинтервалов  той же меры, что у $E_j^i$, получая набор  
$$E_j^i, TE_j^i, T^2 E_j^i,\dots, T^{h_j-1}E_j^i, T^{h_j}E_j^i, T^{h_j+1}E_j^i, \dots, T^{h_j+s_j(i)-1}E_j^i$$
(все эти множества  не пересекаются).
Обозначив $E_{j+1}= E^1_j$, для   $i<r_j$ положим 
$$T^{h_j+s_j(i)}E_j^i = E_j^{i+1}.$$
 Набор надстроеных колонн с этого момента  рассматривается как   башня  этапа $j+1$,  состоящая из полуинтервалов  
$$E_{j+1}, TE_{j+1}, T^2 E_{j+1},\dots, T^{h_{j+1}-1}E_{j+1},$$
где  
 $$ h_{j+1}+1 =(h_j+1)r_j +\sum_{i=1}^{r_j}s_j(i).$$

Частичное определение преобразования $T$ на этапе $j$ сохраняется на всех следующих этапах. В итоге на пространстве  $X=\cup_j X_j$ определено  обратимое преобразование $T:X\to X$, сохраняющее  стандартную меру Лебега на $X$.

Преобразование $T$  и индуцированный им унитарный оператор 
$T$, $Tf(x)=f(Tx)$,  в статье обозначаются одинаково.
Автоморфизм ранга один эргодичен,   имеет простой (однократный) спектр.  Известно, что  индикатор башни $X_1$ являются циклическим вектором для оператора $T$ с непрерывным спектром, что заведомо выполнено, если мера $X$ бесконечна.

\vspace{2mm}
\bf Пример конструкции  $\bf T$, для которой 
 $\bf T^{\otimes 2}$ имеет спектральную кратность 2, 
спектр $\bf T^{\otimes 3}$  лебеговский. 

\rm
 Пусть последовательность $J_k$ определена так:

$J_1=1,\  J_{k+1}= J_k +(k+1)^8.$

\vspace{2mm}
Положим

$r_{j}:=(k+1)^2 $, при $j(k)\leq j<j(k+1)$;

$s_j(1):=2h_j-1$; \ $s_j(2):=2h_j$; 

 $s_j(i):=2^ih_j$ при   $2 <  i \leq r_j,$   

где  $h_1=2$,  $h_{j+1}= r_jh_j +\sum_{i=1}{r_j} s_j(i)$.
 
 Значение  $r_j=(k+1)^2$ повторяется много раз ($r_j^4$), чтобы  для сумм  $$F_k=\sum_{J_k<j< J_{k+1}} {\bf 1}_{T^{h_j}E_1}\otimes {\bf 1}_{T^{h_j}E_1}$$
функции $F_k /\|F_k\|$ асимптотически были близки к индикатору множества 
$$\uu_{0\leq p,q\leq 1}{T^pE_1}\times {T^qE_1}.$$
 Это устанавливается путем   модификации метода работы \cite{R24}. Благодаря указанной аппроксимации  получаем следующее утверждение.
  
\vspace{3mm}
\bf  Лемма. \it  Пусть $f={\bf 1}_{E_1}\times {\bf 1}_{E_1}$.  
Тогда  для всяких $m,n\in \Z$   любое  замкнутое и инвариантное  отосительно оператора $T\otimes T$   пространство вместе с вектором $(T^m\otimes T^n)f$  содержит 
вектор $$(T^{m+1}\otimes T^n)f + (T^m\otimes T^{n+1})f.$$\rm

\vspace{3mm}
Из  леммы  легко вытекает, что  $T\otimes T$-инвариантное пространство, содержащее  векторы $f$ и $(T\otimes I)f$,
содержит все векторы вида $(T^m\otimes T^n)f$. Значит, $T\otimes T$ обладает однородным спектром кратности 2. Следовательно, 
 $T\odot T$ имеет простой спектр, так как  
$$T\otimes T\cong (T\odot T) \,\bigoplus\, (T\odot T).$$

Если $\hat\sigma(n)=\mu(E_1\cap T^nE_1)$, то для коэффциентов 
$$c(n):=\hat\sigma(n)^3=\hat{\sigma^{\ast 3}}(n)$$
непосредственно проверяется сходимость ряда 
$\sum_n c(n)^2$. Следовательно, спектр $T^{\odot 3}$  абсолютно непрерывный, так как спектральный тип оператора $T^{\odot 3}$ есть сверточная степень $\sigma^{\ast 3}$ спектральной меры $\sigma$ оператора $T$. Если сходимость указанного ряда достаточно быстрая, то $\sigma^{\ast 3}$ гарантированно эквивалентна  лебеговской мере. 

Многочисленные  результаты о спектрах сохраняющих меру действияй описаны в недавнем обзоре \cite{KL}.   К ним автор может добавить, например,  конструктивные решения задачи Колмогорова о групповом свойстве спектра и  задачи Рохлина об однородном спектре в классе перемешивающих систем; новые реализации спектров для гауссовских и пуассоновских систем; решение задачи  Оселедца о проекциях тензорного квадрата (спектральной) меры (см. работы 2020-2024 гг).

\vspace{4mm}
\bf  Доказательство теоремы 1. \rm  
Сперва рассмотрим \bf случай $\bf n=2$. \rm На отрезке $[0, 1]$ определим  множества $A,B$:

$$A=\left\{\sum_{i=1}^\infty  \frac {a_i} {4^i} \, :\,  a_i\in\{0,1\}\right\},$$

$$B=2A=\left\{\sum_{i=1}^\infty  \frac {b_i} {4^i} \, :\,  b_i\in\{0,2\}\right\}.$$

 Рассмотрим теперь    меры с носителями   $A$ и $B$, соответственно.
Первая $\sigma_1$ --   бернуллиевская мера типа 
$(\frac 1 2, \frac 1 2,  0,0)$, вторая мера $\sigma_2$ -- бернуллиевская  типа   $(\frac 1 2, 0, \frac 1 2, 0)$.

(i) Отображение $C_2: A\times B\to [0,1]$, 
$$ C_2(a,b)=a+b, $$
является изоморфизмом пространств $([0,1]^2, \sigma_1\times \sigma_2)$
и $([0,1],\lambda)$  ($\lambda$ -- мера Лебега). Это несложно проверить (способ стандартный: цилиндр переходит в  интервал  с сохранением меры).

Пусть $U$ -- оператор умножения, действующий  в пространстве $L_2([0,1], \sigma_1)$, 
$$ Uf(x)=e^{2\pi i x}f(x).$$

Тогда действие оператора $V=U\otimes U^2$  на $([0,1]^2, \sigma_1\times \sigma_1)$ описывается  формулой 
$$ VF(x,x')=e^{2\pi i (x+2x')}F(x,x').$$
Но для почти всякого $y\in [0,1]$ найдутся  единственные $x,x'\in A$ такие, что $x+2x'=y$, поэтому в силу (i) оператор $V$ изоморфен
оператору умножения на $z=e^{2\pi i y}$ в пространстве $L_2([0,1].\lambda)$, что и требовалось доказать.

 \vspace{2mm}
\bf Случай $\bf n=3$.  \rm
Определим 
$$A=\left\{\sum_{i=1}^\infty  \frac {x_i} {8^i} \, :\,  x_i\in\{0,1\}\right\};$$
Следует сказать, что определение множества $A$ зависит от  $n$.

Замечаем, что для почти всякого $y$ найдутся единственные $x,x',x''\in A$
такие, что $y=x+2x'+4x''$, причем  отображение 
$$(x,x',x'')\to x+2x'+4x''$$
является изоморфизмом пространств  
$([0,1]^3, \sigma_1^{\otimes 3})$ и 
и $([0,1],\lambda)$, где 
$\sigma_1$ -- мера Бернулли типа $(\frac 1 2, \frac 1 2, 0,0, 0, 0, 0,0)$. 

Пусть  $U$ теперь обозначает   оператор умножения в пространстве $L_2([0,1], \sigma_1)$:
$$ Uf(x)=e^{2\pi i x}f(x).$$

Тогда оператор $V_3=U\otimes U^2\otimes U^4$ есть умножение 
на  $e^{2\pi i (x +2x'+4x'')}$ в пространстве $([0,1]^3, \sigma_1^{\otimes 3})$, а в силу сказанного он 
имеет простой лебеговский спектр.

 Общий  случай $\bf n>3$ рассматривается аналогично очевидным образом. Теорема доказана.

\vspace{4mm}
\bf Замечания, вопросы. \rm   

 Мы знаем, как растут спектральные кратности степеней $V^n$,
они равны $n$.  Надеюсь, читатель согласится с тем, что любопытно узнать, как  как растут спектральные кратности  степеней $U^n$?
Они не могут асимптотически расти линейно, так как в тензорном произведении мы бы увидели степенной рост.  Аналогичный вопрос
об устройстве  спектра оператора $$Gauss(U)^n=\exp(U^n)= \bigoplus_{d=0}^\infty
 \left(U^{\odot d}\right)^n$$
тоже интересен.

\vspace{2mm}
В \cite{R22} решена задача Оселедца: найдется мера  $\sigma$ такая, что 
плотны направления, когда проекции $\sigma\times\sigma$ абсолютно непрерывны и одновременно плотны направления (на горизонтальныю прямую), когда проекции $\sigma\times\sigma$ сингулярны.

\vspace{2mm}
В.И. Оселедец  высказал интересную гипотезу: \it невозможна ситуация, когда проекция   меры $\sigma\times\sigma$ для некоторого направления  абсолютно непрерывна, а все достаточно близкие направления дают сингулярные меры.\rm
Наш пример, в силу самоподобия меры $\sigma_1$ не дает контрпримера к гипотезе, но, быть может, слегка на него похож.

\vspace{2mm}
Существует ли оператор $U$ такой, что  $U^{\odot n}$ обладает простым абсолютно непрерывным спектром?


\end{document}